\documentclass{article}
\usepackage{latexsym}
\usepackage{amstext}
\usepackage{amsmath}
\usepackage{amssymb}
\usepackage{amsthm}

\newtheorem{theorem}{Theorem}
\newtheorem{lemma}[theorem]{Lemma}
\newtheorem{corollary}[theorem]{Corollary}
\theoremstyle{definition}
\newtheorem{definition}[theorem]{Definition}

\newtheorem*{remark}{Remark}

\def \co {\mathcal{O}}

\def \P {\mathbb{P}}
\def \tf {\tilde{f}}
\def \tg {\tilde{g}}

\begin{document}
\bibliographystyle{amsplain}
\title{One-Parameter Families of Unit Equations}
\author{Aaron Levin\\adlevin@math.brown.edu}
\date{}
\maketitle
\begin{abstract}
We study one-parameter families of $S$-unit equations of the form $f(t)u+g(t)v=h(t)$, where $f$, $g$, and $h$ are univariate polynomials over a number field, $t$ is an $S$-integer, and $u$ and $v$ are $S$-units.  For many possible choices of $f$, $g$, and $h$, we are able to determine all but finitely many solutions to the corresponding one-parameter family of $S$-unit equations.  The results are obtained as consequences of some recent results on integral points on surfaces. 
\end{abstract}
\section{Introduction}
An equation of fundamental interest in number theory is the $S$-unit equation
\begin{equation}
\label{2unit}
au+bv=c \quad \text{ in }u,v\in \co_{k,S}^*,
\end{equation}
where $k$ is a number field, $a,b,c\in k^*$, $S$ is a finite set of places of $k$ containing the archimedean places, $\co_{k,S}$ is the ring of $S$-integers of $k$, and $\co_{k,S}^*$ is the group of $S$-units of $k$.  The basic theorem on the $S$-unit equation (\ref{2unit}) is
\begin{theorem}[Siegel, Mahler]
\label{uth}
The set of solutions to \eqref{2unit} is finite.
\end{theorem}
This was proved by Siegel in the case $S$ consists of the set of archimedean places of $k$ and extended by Mahler to arbitrary $S$.

Equation (\ref{2unit}) and Theorem \ref{uth} have been generalized in at least three distinct directions.  First, it is easy to see that Theorem \ref{uth} is equivalent to the assertion that, in any affine embedding, $\P^1\setminus \{\text{three points}\}$ has only finitely many $S$-integral points.  Thus, Siegel's theorem on integral points on affine curves may be seen as a generalization of Theorem \ref{uth}.  Second, the hypothesis that $u,v\in \co_{k,S}^*$ has been generalized, for instance, by Lang \cite{La} to the assumption that $u,v\in \Gamma$, where $\Gamma$ is a finitely generated subgroup of $\mathbb{C}^*$.  Third, we can consider $S$-unit equations with more variables and terms.  The main theorem in this case, proved independently by Evertse \cite{Ev} and van der Poorten and Schlickewei \cite{vdP2}, is the following.
\begin{theorem}[Evertse, van der Poorten and Schlickewei]
\label{main}
All but finitely many solutions of the equation
\begin{equation*}
\alpha_1u_1+\alpha_2u_2+\ldots+\alpha_nu_n=\alpha_{n+1} \quad \text{ in } u_1,\ldots,u_n\in \co_{k,S}^*,
\end{equation*}
where $\alpha_1,\ldots,\alpha_n\in k^*$, satisfy an equation of the form $\sum_{i\in I}\alpha_iu_i=0$, where $I$ is a subset of $\{0,\ldots,n\}$.
\end{theorem}
Of course, one can also consider combinations of the above generalizations, such as \cite{Ev2} and \cite{vdP}, which extend Theorem \ref{main} to finitely generated subgroups of $\mathbb{C}^*$.  We mention also that Theorem \ref{uth} has been refined in various ways.  There are the effective estimates for the heights of solutions to (\ref{2unit}) coming from linear forms in logarithms \cite{Ev3,Gy} and bounds on the number of solutions to (\ref{2unit}) depending only on the size of $S$ (and in some versions $[k:\mathbb{Q}]$) \cite{Be,Bo,Ev4,Gy}.

In this article we consider another possible generalization of (\ref{2unit}).  We study one-parameter families of the two variable $S$-unit equation, namely,
\begin{equation}
\label{unit}
f(t)u+g(t)v=h(t) \quad \text{ in } t\in \co_{k,S}, \enspace u,v\in \co_{k,S}^*,
\end{equation}
where $f$, $g$, and $h$ are nonzero polynomials in $k[t]$.  In most situations, our requirement that $t\in \co_{k,S}$ can be relaxed to the possibly more natural condition that $t\in k$ (see Lemma \ref{lem}).

Even for very simple choices of $f$, $g$, and $h$, Eq. (\ref{unit}) leads to open problems.  For instance, taking $f=g=1$ and $h=t^2-1$, we obtain the equation
\begin{equation*}
u+v=t^2-1 \quad \text{ in } t\in \co_{k,S}, \enspace u,v\in \co_{k,S}^*.
\end{equation*}
Solving this equation is essentially the same as determining when the sum of three $S$-units is a perfect square.  This appears to be a difficult problem.  Indeed, determining whether or not there are infinitely many perfect squares in $\mathbb{Z}$ of the form $2^a+3^b+1$ for positive integers $a$ and $b$ is already an open problem.  However, as an example of our results we will show 
\begin{theorem}
\label{gen}
For general nonconstant polynomials $f,g,h\in k[t]$ with 
\begin{equation*}
\deg f+\deg g=\deg h>2, 
\end{equation*}
the equation
\begin{equation*}
f(t)u+g(t)v=h(t) \quad \text{ in }t\in k, \enspace u,v\in \co_{k,S}^*,
\end{equation*}
has only finitely solutions with $f(t)g(t)h(t)\neq 0$.
\end{theorem}
By ``general" here, we mean that if one parametrizes the polynomials $f$, $g$, and $h$ in Theorem \ref{gen} in the obvious way by the affine space $\mathbb{A}^{2\deg h+3}$, then we are excluding polynomials $f$, $g$, and $h$ parametrized by some Zariski-closed subset of $\mathbb{A}^{2\deg h+3}$ (in principle, the Zariski-closed subset could be explicitly given).

If $t_0\in \co_{k,S}$ is not a root of $f$, $g$, or $h$, then substituting $t=t_0$ into (\ref{unit}) gives an $S$-unit equation.  Therefore, Theorem \ref{gen} gives numerous examples of $S$-unit equations which have no solution.  However,  since our methods are ineffective, we cannot determine the finitely many values of $t$ for which there are solutions in Theorem \ref{gen}, and so we cannot explicitly determine (by our methods) for any given value of $t$ that the corresponding $S$-unit equation has no solutions.

Our results are proven as consequences of some recent theorems on integral points on surfaces \cite{Co6,Le}.  These theorems trace their origin to the new proof of Siegel's theorem using the Schmidt subspace theorem given in \cite{Co} and developed for surfaces in \cite{Co2}.  The ultimate reliance of our results on the subspace theorem is the reason for their ineffectivity.\\\\
{\bf Acknowledgments.}  I would like to thank Pietro Corvaja and Umberto Zannier, without whom this paper would not have been written.  The ideas involved originated in conversations amongst the three of us while attending the program on Diophantine Geometry at the De Giorgi Center in Pisa, Italy.
\section{Elementary Observations}
Fix nonzero polynomials $f$, $g$, and $h$.  We first discuss an obvious set of solutions to (\ref{unit}).  Let $Z(fgh)$ denote the set of zeroes of $fgh$ and suppose that $S$ is large enough such that $Z(fgh)\subset \co_{k,S}$.  Equation (\ref{unit}) is not a unit equation for $t\in \co_{k,S}$ exactly when $t\in Z(fgh)$.  For these values of $t$, the set of solutions to (\ref{unit}) is easily described.  For instance, if $f(t_0)=0$ and $g(t_0)\neq 0$, then the set of solutions to (\ref{unit}) with $t=t_0$ is given by $u\in \co_{k,S}^*$ and $v=\frac{h(t_0)}{g(t_0)}$, assuming that $\frac{h(t_0)}{g(t_0)}\in \co_{k,S}^*$ (otherwise, there are no solutions with $t=t_0$).  Thus, we will call any solution with $t\in Z(fgh)$ a trivial solution and any solution with $t\notin Z(fgh)$ a nontrivial solution.

In the rest of the paper, we will always make the assumption that $f$ and $g$ do not have a common zero.  We now show that there is no loss of generality in doing this.
\begin{lemma}
Let $f,g,h\in k[t]$ be nonzero polynomials.  There exist polynomials $f',g',h'\in k[t]$ such that $f'$ and $g'$ do not have a common zero and such that, for large enough $S$, there is a natural inclusion of the set of solutions to
\begin{equation}
\label{l1}
f(t)u+g(t)v=h(t) \quad \text{ in }t\in \co_{k,S}, \enspace u,v\in \co_{k,S}^*,
\end{equation}
into the set of solutions to
\begin{equation}
\label{l2}
f'(t')u'+g'(t')v'=h'(t') \quad \text{ in }t'\in \co_{k,S}, \enspace u',v'\in \co_{k,S}^*.
\end{equation}
\end{lemma}
\begin{proof}
We easily reduce to the case where $f$, $g$, and $h$ do not all have a common zero.  Let $d\in k[t]$ be such that $f'=f/d$ and $g'=g/d$ do not have a common zero and $f',g'\in \co_{k,S}[t]$.  For any $(t,u,v)$ satisfying (\ref{l1}) we have $f'(t)u+g'(t)v=h(t)/d(t)\in \co_{k,S}$.  It follows from the fact that $h$ and $d$ do not have a common zero that, after enlarging $S$, $d(t)\in \co_{k,S}^*$ for any $(t,u,v)$ satisfying (\ref{l1}).  Therefore, if $(t,u,v)$ is a solution to (\ref{l1}) then $(t',u',v')=(t,ud(t),vd(t))$ is a solution to (\ref{l2}), where we have set $h'=h$.
\end{proof}

Finally, in most situations the restriction that $t$ is an $S$-integer in (\ref{unit}) is unnecessary.
\begin{lemma}
\label{lem}
Suppose that the largest degree among $f$, $g$, and $h$ is uniquely attained among $f$, $g$, and $h$.  Then for large enough $S$, the set of solutions to
\begin{equation}
\label{unit2}
f(t)u+g(t)v=h(t) \quad \text{ in } t\in k, \enspace u,v\in \co_{k,S}^*,
\end{equation}
is the same as the set of solutions to \eqref{unit}.
\end{lemma}
\begin{proof}
Clearly, we can assume that $f,g,h\in \co_k[t]$.  Suppose that $S$ is large enough such that the leading coefficients of $f$, $g$, and $h$ are $S$-units.  Let $u,v \in \co_{k,S}^*$.  Then by our assumption on the degrees of $f$, $g$, and $h$,  it follows that the leading coefficient of $f(t)u+g(t)v-h(t)$ (as a polynomial in $t$) is an $S$-unit.  Therefore, if $(t,u,v)$ is a solution to (\ref{unit2}) then $t$ must be an $S$-integer.
\end{proof}

\section{Integral Points on Certain Affine Surfaces}
We start with a definition of integral points for affine varieties.
\begin{definition}
Let $V$ be an affine variety defined over a number field $k$.  Let $S$ be a finite set of places of $k$ containing the archimedean places.  We define a set $R\subset V(k)$ to be a set of $S$-integral points on $V$ if there exists an affine embedding $\phi:V\hookrightarrow\mathbb{A}^n$ such that $\phi(R)\subset \mathbb{A}^n(\co_{k,S})$.
\end{definition}
The solutions to (\ref{unit}) are intimately related to integral points on certain affine surfaces.  For our purposes, it will be most convenient to view the affine surfaces of interest as subsets of $\P^1\times \P^1$.
\begin{theorem}
\label{main2}
Let $f,g,h\in k[t]$ be nonzero polynomials and let 
\begin{equation*}
T=\{(t,u,v)\in \co_{k,S}\times \co_{k,S}^*\times \co_{k,S}^*\mid f(t)u+g(t)v=h(t)\}
\end{equation*}
be the set of solutions to \eqref{unit}.  Suppose that $f$ and $g$ do not have a common zero.  Let $\tf,\tg\in k[t]$ be such that $f\tg+g\tf=h$.  Consider $\P^1\times \P^1$ with coordinates $(x_1,y_1)\times (x_2,y_2)$.  Let $Z$ be the closed subset of $\P^1\times \P^1$ that is the union of the sets defined by the four equations (appropriately clearing denominators in the last two equations)
\begin{align*}
y_1&=0,\\
y_2&=0,\\
x_1f\left(\frac{x_2}{y_2}\right)-y_1\tf\left(\frac{x_2}{y_2}\right)&=0,\\
x_1g\left(\frac{x_2}{y_2}\right)-y_1\tg\left(\frac{x_2}{y_2}\right)&=0.
\end{align*}
Let $R\subset \P^1\times\P^1\setminus Z$ be the set
\begin{multline}
\label{Req}
R=\{(\tf(t)-v,f(t))\times (t,1)\mid (t,u,v)\in T, f(t)\neq 0\}\cup \\
\{(u-\tg(t),g(t))\times (t,1)\mid (t,u,v)\in T, g(t)\neq 0\}.
\end{multline}
Then $R$ is a set of $S$-integral points on $\P^1\times \P^1\setminus Z$.
\end{theorem}
\begin{proof}
Multiplying everything by a scalar, we may assume that $f,g,h\in \co_k[t]$.  We first show that there exists a constant $c\in k^*$ such that for all $(t,u,v)\in T$,
\begin{equation}
\label{int}
\frac{c(\tf(t)-v)}{f(t)}\in \co_{k,S} \text{ if } f(t)\neq 0, \quad \frac{c(u-\tg(t))}{g(t)}\in \co_{k,S} \text{ if } g(t)\neq 0.
\end{equation}
By possibly making $c$ larger, it clearly suffices to prove this for all but finitely many values of $t$.  So we will ignore values of $t$ for which $f(t)=0$ or $g(t)=0$.  In this case, it follows from (\ref{unit}) and the definitions of $\tf$ and $\tg$ that 
\begin{equation}
\label{eq1}
\frac{\tf(t)-v}{f(t)}=\frac{u-\tg(t)}{g(t)}.
\end{equation}
Since $f$ and $g$ do not have a common zero, there exist polynomials $p_1,p_2\in \co_k[t]$ such that $fp_1+gp_2=a$, where $a\in \co_k$ is a constant.  Let $b_1,b_2\in \co_k$ be such that $b_1\tf$ and $b_2\tg$ have integral coefficients.  Then it follows from (\ref{eq1}) and the fact that $t,u,v\in \co_{k,S}$ that we can take $c=ab_1b_2$ in (\ref{int}) if $f(t)g(t)\neq 0$.

Since $y_1y_2\neq 0$ on $V=\P^1\times\P^1\setminus Z$, let $x_1'=x_1/y_1$ and $x_2'=x_2/y_2$ be coordinates on $V$.  Then every regular function on $V$ may we be written as $p(x_1',x_2')/((x_1'f(x_2')-\tf(x_2'))^m(x_1'g(x_2')-\tg(x_2'))^n)$, where $p$ is a polynomial in two variables and $m$ and $n$ are integers.  A simple calculation shows that for $(t,u,v)\in T$, if $x_1'=(\tf(t)-v)/f(t)$ or $x_1'=(u-\tg(t))/g(t)$, and $x_2'=t$, then 
\begin{align*}
&x_1'f(x_2')-\tf(x_2')=-v,\\
&x_1'g(x_2')-\tg(x_2')=u.
\end{align*}
For these values of $x_1'$ and $x_2'$, it follows from (\ref{int}) that there exists a constant $d\in k^*$ such that $dp(x_1',x_2')\in \co_{k,S}$.  Therefore, for any regular function $\psi$ on $V$ we see that there exists a constant $d\in k^*$ such that $d\psi(R)\subset \co_{k,S}$.  Note also that $V$ is affine.  Thus, after multiplying the coordinate functions by suitable constants, for any affine embedding $\phi:V\hookrightarrow \mathbb{A}^N$ we have $\phi(R)\subset \mathbb{A}^N(\co_{k,S})$.
\end{proof}
So the problem of determining solutions to (\ref{unit}) is now reduced to the study of $S$-integral points on certain affine surfaces.  When there does not exist a Zariski-dense set of $S$-integral points on such a surface, we can parametrize the solutions to the corresponding one-parameter $S$-unit equation.
\begin{theorem}
Let $f$, $g$, $h$ and $Z\subset \P^1\times \P^1$ be as in Theorem \ref{main2}.  Suppose that there does not exist a Zariski-dense set of $S$-integral points on $\P^1\times \P^1\setminus Z$.  Then there exist finitely many quintuples $(z_i,a_i,b_i,p_i,q_i)$, $z_i\in k[t,\frac{1}{t}]$, $a_i,b_i\in k$, $p_i,q_i \in \mathbb{Z}$, with
\begin{equation}
\label{zeq}
a_if(z_i(t))t^{p_i}+b_ig(z_i(t))t^{q_i}=h(z(t))
\end{equation}
for $i=1,\ldots,j$, such that all solutions to \eqref{unit} are parametrized by
\begin{equation}
\label{par}
t=z_i(s), u=a_is^{p_i}, v=b_is^{q_i}, s\in k
\end{equation}
for $i=1,\ldots, j$.
\end{theorem}
This follows easily from Siegel's theorem.
\begin{proof}
Let $R$ be as in (\ref{Req}).  Then by Theorem \ref{main}, $R$ is a set of $S$-integral points on $\P^1\times \P^1\setminus Z$.  By hypothesis, $R$ is not Zariski-dense.  Let $C_i$, $i=1,\ldots, j$ be the one-dimensional irreducible components of the Zariski-closure of $R$ in $\P^1\times \P^1$.  By Siegel's theorem, $C_i$ is a rational curve defined over $k$, and if $\phi_i:\P^1\to C_i\to \P^1\times \P^1$ is the normalization map composed with the inclusion map of $C_i$, then $\#\phi_i^{-1}(Z\cap C_i)\leq 2$.  After an automorphism of $\P^1$, we can assume that $\phi_i^{-1}(Z\cap C_i)\subset \{0,\infty\}\subset \P^1$.  Let $\phi_i(t)=(y_i(t),1)\times (z_i(t),1)$ in affine coordinates on $\P^1$.  By the definition of $Z$ and our assumption that $\phi_i^{-1}(Z\cap C_i)\subset \{0,\infty\}$ we then have $y_i,z_i\in k[t,\frac{1}{t}]$ and
\begin{align*}
y_i(t)f(z_i(t))-\tf(z_i(t))&=a_it^{p_i},\\
y_i(t)g(z_i(t))-\tg(z_i(t))&=b_it^{q_i}
\end{align*}
for some $a_i,b_i\in k$ and some integers $p_i,q_i\in \mathbb{Z}$.  Now easy calculations and the definition of $R$ show that (\ref{zeq}) holds and that all but finitely many solutions to (\ref{unit}) are parametrized by (\ref{par}) for $i=1,\ldots,j$.  The finitely many remaining solutions can be covered in the theorem by taking, for some $i'$,  $z_{i'}(s)$ constant and $p_{i'}=q_{i'}=0$ with appropriate $a_{i'}, b_{i'}\in k$.
\end{proof}
We define a set of curves on a surface to be in general position if the intersection of any three of the curves is empty.  Recall also that a curve $C$ on $\P^1\times \P^1$ is said to be of type $(a,b)$ if it is defined by a bihomogeneous equation of bidegree $(a,b)$.  We will need the following theorem on integral points from \cite{Le}.
\begin{theorem}
\label{myth}
Let $Z_1,Z_2,Z_3,Z_4\subset \P^1\times \P^1$ be curves in general position of types $(0,1)$, $(1,0)$, $(1,m)$, and $(1,n)$, respectively.  Then there exists a Zariski-closed subset $Y\subset \P^1\times \P^1$, independent of $k$ and $S$, such that for any set $R$ of $S$-integral points on $\P^1\times \P^1\setminus \cup_{i=1}^4Z_i$ the set $R\setminus Y$ is finite.
\end{theorem}
More generally, the theorem holds with $Z_3$ and $Z_4$ of types $(a,b)$ and $(c,d)$, respectively, with $a,b,c,d>0$.
We now give some more information on the exceptional set $Y$ in Theorem \ref{myth}.  It is easy to see that every curve intersects $Z=\cup_{i=1}^4Z_i$ in at least two points.  So by Siegel's theorem, every irreducible curve $C$ in (a minimal) $Y$ intersects $Z$ in exactly two points $P$ and $Q$.  We denote the intersection number of two curves $D$ and $E$ on a surface by $D.E$.
\begin{theorem}
\label{exc}
Suppose $m\geq n$.  The set $Y$ consists of the following types of irreducible curves $C$:
\begin{enumerate}
\item $C$ is of type $(0,1)$ or $(1,0)$.
\item $C$ is of type $(1,p)$ or $(q,1)$ with $p,q>0$, $p\leq m$, $(q-1)n\leq m-1$, and $P\in Z_1\cap Z_2$, $Q\in Z_3\cap Z_4$ (up to switching $P$ and $Q$).
\label{2}
\item $C$ is of type $(1,p)$ with $0<p\leq m$ and either $P\in Z_1\cap Z_3, Q\in Z_2\cap Z_4$ or $P\in Z_1\cap Z_4, Q\in Z_2\cap Z_3$ (up to switching $P$ and $Q$).
\label{3}
\end{enumerate}
Furthermore, if $m+n>2$, then for general $Z_3$ and $Z_4$ of types $(1,m)$ and $(1,n)$, respectively, $Y$ consists only of $(0,1)$ curves.
\end{theorem}
\begin{proof}
By Siegel's theorem, $C$ must be a rational curve and nonsingular at $P$ and $Q$.  Suppose $C$ is of type $(a,b)$ with $a,b>0$.  Then $C$ intersects each of $Z_1$, $Z_2$, $Z_3$, and $Z_4$ in at least one point.  Therefore, by the general position assumption, $P$ and $Q$ must be as in \ref{2} or \ref{3}.  Suppose $P$ and $Q$ are as in \ref{2}.  Since $Z_1$ and $Z_2$ intersect transversally at $P$, $C$ is nonsingular at $P$, and $C$ intersects $Z_1$ and $Z_2$ in exactly one point, we must have either $C.Z_1=1$ or $C.Z_2=1$.  Therefore $C$ is of type $(1,p)$ or $(q,1)$ with $p,q>0$.  To prove the inequalities on $p$ and $q$, we use the intersection formula \cite[Ch. V:Ex. 3.2]{Ha}
\begin{equation}
\label{inter}
D.E=\sum \mu_{P'}(D)\mu_{P'}(E)
\end{equation}
for curves $D$ and $E$, where $\mu_{P'}(D)$ and $\mu_{P'}(E)$ denote the multiplicity of the point $P'$ on $D$ and $E$, respectively, and the sum is taken over all infinitely near points $P'$ on $\P^1\times \P^1$.  Assume that $Z_3$ and $Z_4$ are irreducible and hence nonsingular (the reducible case is similar).  Let $Q_1,Q_2,\ldots$ be the infinitely near points of $C$ infinitely near $Q$ where $Q_1$ lies on the blow-up $X_1$ of $\P^1\times \P^1$ at $Q$, $Q_2$ lies on the blow-up $X_2$ of $X_1$ at $Q_1$, and so on.  Then it follows from (\ref{inter}) that $\mu_{Q_i}(Z_3)=1$ for $i=1,\ldots,C.Z_3-1$ and $\mu_{Q_i}(Z_4)=1$ for $i=1,\ldots,C.Z_4-1$.  Therefore, again from (\ref{inter}), we obtain $Z_3.Z_4\geq \min\{C.Z_3,C.Z_4\}$.  This is equivalent to $p\leq m$ and $(q-1)n\leq m-1$.

Now suppose that $P\in Z_1\cap Z_3$ and $Q\in Z_2\cap Z_4$.  Since $Z_1$ intersects $Z_3$ transversally at $P$, $C.Z_3>1$, and $C$ intersects $Z_1$ and $Z_3$ only at $P$, it follows that $C.Z_1=1$, i.e., $C$ is of type $(1,p)$.  Using (\ref{inter}) at the point $Q$ as before, we obtain $p\leq m$.  The other case for $P$ and $Q$ follows similarly.

Now suppose that $m+n>2$.  It is not hard to see that for $Y$ to contain curves other than $(0,1)$ curves, the curves in $Z$ must be in special position, definable by algebraic relations.  For example, if $Z_1$, $Z_2$, $Z_3$ and $Z_4$ meet pairwise transversally, then case \ref{2} of the theorem can never occur.  We leave the details to the reader.
\end{proof}
As a consequence of Theorems \ref{myth} and \ref{exc} we obtain
\begin{corollary}
\label{cor1}
Let $f,g,h\in k[t]$ be nonconstant polynomials such that $f$ and $g$ do not have a common zero and $\deg f+\deg g=\deg h$.  Let $m=\deg f$ and $n=\deg g$ and suppose that $m\geq n$.  Then all but finitely many solutions to
\begin{equation}
\label{unit4}
f(t)u+g(t)v=h(t) \quad \text{ in } t\in k, \enspace u,v\in \co_{k,S}^*,
\end{equation}
are parametrized by a finite number of families, independent of $k$ and $S$, of the form
\begin{equation*}
t=z(s), u=as^p, v=bs^q, \quad s\in k
\end{equation*}
where $z\in k[t]$, $a,b\in k$, $p,q\in \mathbb{Z}$,
\begin{equation*}
af(z(t))t^p+bg(z(t))t^q=h(z(t)),
\end{equation*}
and $(\deg z-1)n\leq m-1$.  Furthermore, $p,q>0$ if $\deg z>1$.
\end{corollary}
\begin{proof}
First note that there exist $\tf, \tg\in k[t]$ with $\deg \tf\leq n$ and $\deg \tg\leq m$ such that $f\tg+g\tf=h$.  To see this, let $P_i$ denote the vector space of polynomials over $k$ of degree at most $i$ and consider the map $P_n\oplus P_m\to P_{m+n}$ given by $x\oplus y\mapsto fx+gy$.  The kernel is one-dimensional, generated by $(-g)\oplus f$ (since $f$ and $g$ do not have a common zero), and so by counting dimensions we see that the map is surjective.  Therefore, by Theorem \ref{main2} with this $\tf$ and $\tg$, we see that solutions to (\ref{unit4}) give rise to a set of $S$-integral points on $\P^1\times\P^1\setminus Z$, where $Z=\cup_{i=1}^4Z_i$ and $Z_1$, $Z_2$, $Z_3$, and $Z_4$ are of types $(1,0)$, $(0,1)$, $(1,m)$, and $(1,n)$, respectively.  Furthermore, it is easy to see that $\deg f+\deg g=\deg h$ implies that $Z_1$, $Z_2$, $Z_3$, and $Z_4$ are in general position (in fact, the other direction also holds).  Now a straightforward translation of Theorems \ref{myth} and \ref{exc} into information about (\ref{unit4}), via the correspondence in Theorem~\ref{main2}, gives the corollary.
\end{proof}
Theorem \ref{gen} from the introduction is similarly a direct consequence of the last statement of Theorem \ref{exc}.  We merely note that a $(0,1)$ curve in $Y$ corresponds to a trivial set of solutions to (\ref{unit}).  As an example of Corollary \ref{cor1}, we explicitly work out what happens when $f$ and $g$ are linear and $h$ is quadratic.  In this case the $z$ in Corollary \ref{cor1} must be linear.  The calculations are then straightforward.
\renewcommand{\labelenumi}{(\alph{enumi}).}
\begin{corollary}
Let $L_1=a_1t+a_0$ and $L_2=b_1t+b_0$ be linear over $k$ with $L_1/L_2$ nonconstant.  Let $Q=c_2t^2+c_1t+c_0\in k[t]$ be quadratic.  Consider the equation
\begin{equation}
\label{unit3}
L_1(t)u+L_2(t)v=Q(t) \quad \text{ in } t\in k, \enspace u,v\in \co_{k,S}^*.
\end{equation}
Let $r_1$ and $r_2$ be the roots of $Q$.  Then there exist the following four families of (potential) solutions to \eqref{unit3}:
\begin{align}
\label{eq2}
&t=\frac{(a_1b_0-a_0b_1)\eta}{c_2(b_1r_1+b_0)}+r_2,\quad u=\eta,\quad v=-\frac{(a_1r_1+a_0)\eta}{b_1r_1+b_0},\quad \eta\in \co_{k,S}^*\\
\label{eq3}
&t=\frac{(a_1b_0-a_0b_1)\eta}{c_2(b_1r_2+b_0)}+r_1,\quad u=\eta,\quad v=-\frac{(a_1r_2+a_0)\eta}{b_1r_2+b_0},\quad \eta\in \co_{k,S}^*\\
\label{eq4}
&t=\frac{a_1\eta}{c_2}+\frac{a_1b_1c_0-a_1b_0c_1+a_0b_0c_2}{c_2(a_1b_0-a_0b_1)}, u=\eta, v=\frac{a_1^2c_0-a_0a_1c_1+a_0^2c_2}{a_1(a_1b_0-a_0b_1)}, \eta\in \co_{k,S}^*\\
\label{eq5}
&t=\frac{b_1\eta}{c_2}+\frac{a_1b_1c_0-a_0b_1c_1+a_0b_0c_2}{c_2(a_0b_1-a_1b_0)}, u=\frac{b_1^2c_0-b_0b_1c_1+b_0^2c_2}{b_1(a_0b_1-a_1b_0)}, v=\eta, \eta\in \co_{k,S}^*
\end{align}
All but finitely many solutions to \eqref{unit3} are given as follows:
\begin{enumerate}
\item If $Q$ is not a perfect square and not of the form $\alpha L_1L_2+\beta, \alpha,\beta\in k$, then all but finitely many nontrivial solutions to \eqref{unit3} are contained in \eqref{eq2}--\eqref{eq5}.
\item If $Q$ is a perfect square with double root $r=r_1=r_2$, then all but finitely many nontrivial solutions to \eqref{unit3} are contained in \eqref{eq2}--\eqref{eq5} and the following family:
\begin{equation*}
t=\eta\sqrt{\frac{a_0b_1-a_1b_0}{b_1c_2}}+r,\quad u=\eta^2,\quad v=-\frac{a_1\eta^2}{b_1}, \quad\eta\in \co_{k,S}^*
\end{equation*}
\item If $Q=\alpha L_1L_2+\beta, \alpha,\beta\in k$, then all but finitely many nontrivial solutions to \eqref{unit3} are contained in \eqref{eq2}--\eqref{eq5} and the following two families:
\begin{align*}
&t=\frac{a_1\eta}{c_2}-\frac{b_0}{b_1},\quad u=\eta,\quad v=\frac{c_2(a_1b_1c_0-a_0b_0c_2)}{a_1^2b_1^2\eta},\quad \eta\in \co_{k,S}^*\\
&t=\frac{b_1\eta}{c_2}-\frac{a_0}{a_1},\quad u=\frac{c_2(a_1b_1c_0-a_0b_0c_2)}{a_1^2b_1^2\eta},\quad v=\eta,\quad \eta\in \co_{k,S}^*
\end{align*}
\end{enumerate}
\end{corollary} 
We can also prove a result in one case where $\deg f+\deg g\neq \deg h$.  We need the following special case of a result from \cite{Co6}.
\begin{theorem}[Corvaja, Zannier]
\label{CZ}
Let $Z_1,Z_2,Z_3,Z_4\subset \P^1\times \P^1$ be curves of types $(0,1)$, $(1,0)$, $(1,1)$, and $(1,1)$, respectively.  Suppose that there exists a unique point where $Z_1$, $Z_3$, and $Z_4$ intersect transversally and that outside of this point of triple intersection the $Z_i$ are in general position.  Then there exists a Zariski-closed subset $Y\subset \P^1\times \P^1$, independent of $k$ and $S$, such that for any set $R$ of $S$-integral points on $\P^1\times \P^1\setminus \cup_{i=1}^4Z_i$ the set $R\setminus Y$ is finite.
\end{theorem}
\begin{remark}
The assumptions of this theorem do not quite satisfy the assumptions of Corollary 1 in \cite{Co6} (with $D_1=Z_1\cup Z_2$, $D_2=Z_3$, and $D_3=Z_4$ in their notation).  However, the proof in \cite{Co6} shows that instead of assuming $D_1$, $D_2$, and $D_3$ are irreducible, it is sufficient that the strict transforms of $D_1$, $D_2$, and $D_3$ in the blow-up at the point of triple-intersection be linearly equivalent to irreducible effective divisors, which certainly occurs in our situation.
\end{remark}
Let $P$ be the point of triple intersection in Theorem \ref{CZ}.  It is easily seen that the $Y$ in Theorem \ref{CZ} can be taken to consist of curves of type $(1,0)$ and $(0,1)$, a curve of type $(1,1)$ tangent to $Z_3$ and passing through $P$ and $Z_2\cap Z_4$, and a curve of type $(1,1)$ tangent to $Z_4$ and passing through $P$ and $Z_2\cap Z_3$.  Using Theorem \ref{main2} to translate this into arithmetic, we obtain
\begin{corollary}
\label{linear}
Let $L_1=a_1t+a_0$, $L_2=b_1t+b_0$, and $L_3=c_1t+c_0$ be linear over $k$ with $L_1/L_2$ nonconstant.  All but finitely many nontrivial solutions to
\begin{equation*}
L_1(t)u+L_2(t)v=L_3(t) \quad \text{ in } t\in \co_{k,S}, \enspace u,v\in \co_{k,S}^*.
\end{equation*}
are parametrized by the following four families:
\begin{align*}
&t=\frac{(a_0b_1-a_1b_0)\eta}{b_1c_1}-\frac{c_0}{c_1}, &&u=\eta, &&v=-\frac{a_1\eta}{b_1}, &&\eta\in \co_{k,S}^*\\
&t=\frac{a_1c_0-a_0c_1-a_1b_0\eta}{a_1b_1\eta}, &&u=\frac{c_1}{a_1}, &&v=\eta,  &&\eta\in \co_{k,S}^*\\
&t=\frac{b_1c_0-b_0c_1-a_0b_1\eta}{a_1b_1\eta}, &&u=\eta, &&v=\frac{c_1}{b_1},  &&\eta\in \co_{k,S}^*\\
&t\in \co_{k,S}, &&u=\frac{b_0c_1-b_1c_0}{a_1b_0-a_0b_1}, &&v=\frac{a_0c_1-a_1c_0}{a_0b_1-a_1b_0}
\end{align*}
\end{corollary}
We note that Corollary \ref{linear} is also implicit in Theorem 2 of \cite{Co6}.
\bibliography{unit}
\end{document}